\documentstyle[11pt]{article}
 \newlength{\standardunitlength}
\newtheorem{cor}{Corollary} \newtheorem{lemma}{Lemma}
\newtheorem{theorem}{Theorem} 
\newenvironment{proof}{\noindent {\sc Proof:}}{$\Box$ \vspace{2 ex}}

\begin{document}

\begin{center} A New Bound for Kloosterman Sums
\end{center}

\begin{center}
By Jason Fulman
\end{center}

\begin{center}
Version of June 20, 2001
\end{center}

\begin{center}
Stanford University (until 8/1/01)
\end{center}

\begin{center}
US Airways (on 8/1/01)
\end{center}

\begin{center}
[University of Pittsburgh (after 8/1/01)]
\end{center}

\begin{center}
Department of Mathematics
\end{center}

\begin{center}
Building 380, MC 2125
\end{center}

\begin{center}
Stanford, CA 94305-2125
\end{center}

\begin{center}
email:fulman@math.stanford.edu
\end{center}

\newpage \begin{abstract} We give generating functions for Gauss sums
for finite general linear and unitary groups. For the general linear
case only our method of proof is new, but we deduce a bound on
Kloosterman sums which is sometimes sharper than Deligne's bound from
algebraic geometry. \end{abstract}

\section{Introduction}

	The problem of bounding exponential sums such as the
Kloosterman sum \[ K(c,d)=\sum_{x \in F_p^*} e^{\frac{2 \pi i}{p}
(cx+ \frac{d}{x})} \] is mathematically central, with applications
to coding theory \cite{H}, spectral graph theory \cite{T}, and modular
forms \cite{Sa}. A main method for bounding such sums is to use deep
results from algebraic geometry.

	Let $\lambda$ be a nontrivial additive character for the
finite field $F_q$, and let $F_q^*$ denote the non-zero elements of
$F_q$. Deligne \cite{D} (see also the exposition \cite{Se}) proves
that for $x \in F_q^*$ \[ \left| \sum_{\alpha_1,\cdots,\alpha_{n-1}
\in F_q^*} \lambda(\alpha_1+\cdots+\alpha_{n-1}+\frac{x}{\alpha_1
\cdots \alpha_{n-1}}) \right| \leq \frac{n}{\sqrt{q}} q^{n/2} .\] As
Section \ref{GaussGL} indicates, it is a simple consequence of Fourier
analysis that \[\left| \sum_{\alpha_1,\cdots,\alpha_{n-1} \in F_q^*}
\lambda(\alpha_1+\cdots+\alpha_{n-1}+\frac{x}{\alpha_1 \cdots
\alpha_{n-1}}) \right| \leq (1-\frac{1}{q-1}) q^{n/2} +
\frac{1}{q-1},\] which is sometimes a stronger bound (though useless
for the case $n=2$).

	The main point of this note is to discuss the relationship of
the second bound with the finite classical groups. Section
\ref{GaussGL} derives a formula for the exponential sum \[ \sum_{g \in
GL(n,q)} \chi(det (g)) \lambda( tr (g)) \] where $\chi$ is a
multiplicative character of $F_q^*$, $\lambda$ is an additive
character of $F_q$, $det$ denotes determinant and $tr$ denotes
trace. This result is known (\cite{E},\cite{Ki1}, \cite{Ko},\cite{L}),
but the proof given here uses the cycle index generating function of
$GL(n,q)$ which has the attractive point of involving ``Euler
products'' over irreducible polynomials.

	Section \ref{GaussGL} then derives the second bound on
Kloosterman sums two paragraphs back by studying characteristic
polynomials of random elements of $SL(n,q)$. Although the derivation
of the bound does not require the study of characteristic polynomials
of random matrices over finite fields, such a derivation is of
interest given that Deligne's work is related to characteristic
polynomials of random matrices from compact Lie groups
\cite{Ka},\cite{KaS}. It would be very interesting to go directly from
finite classical groups to compact Lie groups; p-adic groups may give
the bridge.

	We note that the link between Kloosterman sums and finite
classical groups underlying the second bound is due to Kim \cite{Ki1},
who used the Bruhat decomposition of $SL(n,q)$ with respect to a
maximal parabolic subgroup. Thus Section \ref{GaussGL} gives a
different explanation of this link, using conjugacy classes (and
Fourier analysis) rather than a Bruhat decomposition. Section
\ref{GaussGL} closes by applying our method to find a simple
generating function for the Kloosterman sum \[ \sum_{g \in GL(n,q)}
\lambda(x tr(g) + y tr(g^{-1}))\] which was evaluated by Kim
\cite{Ki4}. Section \ref{GaussU} derives a generating function for the
exponential sum \[ \sum_{g \in U(n,q) \subset GL(n,q^2)} \chi(det (g))
\lambda( tr (g)) \] where $\chi$ is a multiplicative character of
$F_{q^2}^*$, and $\lambda$ is a nontrivial additive character of
$F_{q^2}$. Kim \cite{Ki2},\cite{Ki3} found quite different formulas
for this exponential sum as involved sums of Kloosterman sums.

	As a final remark, our motivation for trying to bound
Kloosterman sums by studying characteristic polynomials of random
matrices arose from card shuffling \cite{F3},\cite{F4}. Two methods of
shuffling, ``affine shuffles'' and ``shuffles followed by cuts'' lead
to exceptionally close distributions on permutations, exactly equal up
to lower order terms involving Ramanujan sums (a type of exponential
sum). The cycle type of a permutation after a $q$ affine shuffle on
$S_n$ has the same distribution as the cycle type of a random degree
$n$ monic polynomial over $F_q$ with constant term 1 and the cycle
type of a permutation after a $q$ riffle shuffle followed by a cut has
the same distribution as the cycle type of a random degree $n$
polynomial over $F_q$ with non-zero constant term. It would be
interesting to find physically natural shuffles whose distribution on
cycle types agrees with that arising from characteristic polynomials
of random matrices.

\section{Gauss Sums for $GL(n,q)$} \label{GaussGL}

	Let $\chi$ be a multiplicative character of $F_q^*$ and let
$\lambda$ be an additive character of $F_q$. Given a polynomial
$P=x^n+c_{n-1}x^{n-1}+\cdots+c_0$ over $F_q$ we define
$\chi(P)=\chi((-1)^n c_0)$ for $c_0 \neq 0$ and $\chi(P)=0$ if
$c_0=0$. We define $\lambda(P)$ as $\lambda(-c_{n-1})$. Note that
$\chi(P_1P_2)=\chi(P_1) \chi(P_2)$ and that
$\lambda(P_1P_2)=\lambda(P_1) \lambda(P_2)$. Also $G(\chi,\lambda)$
denotes the Gauss sum

\[ \sum_{x \in F_q^*} \chi(x) \lambda(x).\]
	
	Before proceeding we recall an elementary and well-known
lemma. In what follows $\phi(z)$ denotes a monic non-constant
irreducible polynomial over $F_q$. For completeness we include a proof
for the case that $\lambda$ is non-trivial.

\begin{lemma} \label{Lfunc}

\[ \prod_{\phi \neq z} \left(\frac{1}{1-\chi(\phi)
\lambda(\phi) u^{deg(\phi)}/q^{i \cdot deg(\phi)}} \right) = \left\{ \begin{array}{ll} \frac{1-u/q^i}{1-u/q^{i-1}} & \mbox{if \ $\lambda,\chi$ \ trivial} \\	1 & \mbox{if \ $\lambda$ \ trivial, $\chi$ \ nontrivial}\\
 1+u G(\chi,\lambda)/q^i & \mbox{if \  $\lambda$ \ nontrivial}
\end{array} \right. \]

\end{lemma} 

\begin{proof} Let $P$ denote a monic polynomial (not necessarily
irreducible) over $F_q$. The left hand side is equal to

\[ \sum_{P: P(0) \neq 0} \frac{u^{deg(P)} \chi(P) \lambda(P)}{q^{i
\cdot deg(P)}}.\] The coefficient of $u^n$ ($n \geq 2$) vanishes
because given $c_0$ the distribution of $c_{n-1}$ for a random monic
degree $n$ polynomial is uniform over all $q$ possible
values. \end{proof}
	
	The first part of Theorem \ref{expGL} is known
(\cite{E},\cite{Ko}\cite{Ki1},\cite{L}) but the method of proof we
give is new. The second part of Theorem \ref{expGL} is simply a
generating function version of the first part.

\begin{theorem} \label{expGL}
\begin{enumerate}

\item \[  1 + \sum_{n=1}^{\infty} \frac{u^n}{|GL(n,q)|} \sum_{g \in
GL(n,q)} \chi(det(g)) \lambda(tr(g)) = \left\{  \begin{array}{ll} \frac{1}{1-u} & \mbox{if \ $\lambda,\chi$ \ trivial}\\
1 & \mbox{$\lambda$ \ trivial, $\chi$ \ nontriv.}\\
\prod_{i \geq 1} (1+\frac{u G(\chi,\lambda)}{q^i}) & \mbox{if  \  $\lambda$ \ nontrivial}.
\end{array} \right.\]

\item For $n \geq 1$, \[ \sum_{g \in GL(n,q)} \chi(det(g)) \lambda(tr(g)) = \left\{ \begin{array}{ll} |GL(n,q)| & \mbox{if \ $\lambda,\chi$ \ trivial}\\
0 & \mbox{if \  $\lambda$ \ trivial, $\chi$ \ nontrivial}\\
q^{n \choose 2} G(\chi,\lambda)^n & \mbox{if  \  $\lambda$ \ nontrivial}.
\end{array} \right. \]
\end{enumerate}
\end{theorem}

\begin{proof} To prove Theorem \ref{expGL} we recall some work of
Stong \cite{St} on the cycle index of $GL(n,q)$ (a survey of
applications of cycle indices of the finite classical groups can be
found in \cite{F2}). As the textbook \cite{He} explains in the section
on rational canonical forms, the conjugacy classes of $GL(n,q)$
correspond to the following combinatorial data: to each monic
non-constant irreducible polynomial $\phi$ over $F_q$, associate a
partition (perhaps the trivial partition) $\nu_{\phi}$ of a
non-negative integer $|\nu_{\phi}|$. We write $\nu \vdash j$ if $\nu$
is a partition of the integer $j$. The only restrictions necessary for
this data to represent a conjugacy class are

\begin{enumerate}
\item $|\nu_z| = 0$
\item $\sum_{\phi} |\nu_{\phi}| deg(\phi) = n$.
\end{enumerate}

	The size of the conjugacy class corresponding to the data
$\nu_{\phi}$ is equal to $\frac{|GL(n,q)|}{\prod_{\phi}
c_{\phi}(\nu_{\phi})}$ where $c_{\phi}(\nu_{\phi})$ is a function of
$\nu$ and $\phi$ which depends on $\phi$ only through the degree of
$\phi$. Define $a_{\phi,\nu}(g)$ to be one if $\nu$ is the partition
corresponding to $\phi$ in the rational canonical form of $g$ and to
be zero otherwise.

	It follows that

\begin{eqnarray*}
& &  1 + \sum_{n=1}^{\infty} \frac{u^n}{|GL(n,q)|} \sum_{g \in
GL(n,q)} \chi(det(g)) \lambda(tr(g)) \\
& = & \prod_{\phi \neq z} \left(1 + \sum_{j=1}^{\infty} \sum_{\nu
\vdash j} \frac{(\chi(\phi) \lambda(\phi) u^{deg(\phi)})^j}{c_{\phi}(\nu_{\phi})} \right).
\end{eqnarray*}

	We now use the fact (derived in \cite{St}) that

\[ 1+\sum_{j=1}^{\infty} \sum_{\nu \vdash j} \frac{(u) ^{j \cdot
deg(\phi)}}{c_{\phi}(\nu_{\phi})} = \prod_{i \geq 1}
\left(\frac{1}{1-u^{deg(\phi)}/q^{i \cdot deg(\phi)}} \right) .\]
Other elementary derivations can be found in \cite{F2}. (There are
similar factorizations for all irreducible polynomials for all finite
classical groups).

	Consequently

\begin{eqnarray*}
& &  1 + \sum_{n=1}^{\infty} \frac{u^n}{|GL(n,q)|} \sum_{g \in
GL(n,q)} \chi(det(g)) \lambda(tr(g)) \\
& = & \prod_{\phi \neq z} \prod_{i \geq 1} \left(\frac{1}{1-\chi(\phi) \lambda(\phi) u^{deg(\phi)}/q^{i \cdot deg(\phi)}} \right) \\
& = &  \prod_{i \geq 1}  \prod_{\phi \neq z}  \left(\frac{1}{1-\chi(\phi) \lambda(\phi) u^{deg(\phi)}/q^{i \cdot deg(\phi)}}\right). \\
\end{eqnarray*}

	Part 1 of Theorem \ref{expGL} now follows from Lemma
\ref{Lfunc}. Part 2 follows from Part 1 using Euler's identity \[
\prod_{i \geq 1} (1+\frac{u}{q^i}) = \sum_{n=0}^{\infty}
\frac{u^n}{(q^n-1) \cdots (q-1)}.\] \end{proof}

	Next we consider the Kloosterman sum

\[ f_{\lambda}(x)=\sum_{\alpha_1,\cdots,\alpha_{n-1} \in F_q^*}
\lambda(\alpha_1+\cdots+\alpha_{n-1}+\frac{x}{\alpha_1 \cdots
\alpha_{n-1}}) \] where $x \in F_q^*$. From pages 46-7 of \cite{Ka},
the Fourier transform of the function $f_{\lambda}$ at the
multiplicative character $\chi$ is equal to
$G(\lambda,\chi)^n$. Fourier inversion implies that

\[ f_{\lambda}(x)= \frac{1}{q-1} \sum_{\chi} \bar{\chi}(x)
G(\lambda,\chi)^n.\] Since $|G(\lambda,\chi)|=q^{n/2}$ for
$\lambda,\chi$ non-trivial \cite{LN}, it follows that \[\left|
\sum_{\alpha_1,\cdots,\alpha_{n-1} \in F_q^*}
\lambda(\alpha_1+\cdots+\alpha_{n-1}+\frac{x}{\alpha_1 \cdots
\alpha_{n-1}}) \right| \leq (1-\frac{1}{q-1}) q^{n/2} +
\frac{1}{q-1}.\] This proves the bound stated in the introduction;
furthermore the technique clearly works for the other sums listed on
page 47 of \cite{Ka}.

	Corollary \ref{GLway} derives the same bound by studying the
trace of characteristic polynomials of elements of $GL(n,q)$ with
determinant $x$. (The first version of this paper only stated the
bound of Corollary \ref{GLway} for $x=1$; Robin Chapman asked us if
the proof could be extended for $x \neq 1$).

\begin{cor} \label{GLway} Let $\lambda$ be a non-trivial additive
character of $F_q$. Then for $x \in F_q^*$, \[\left|
\sum_{\alpha_1,\cdots,\alpha_{n-1} \in F_q^*}
\lambda(\alpha_1+\cdots+\alpha_{n-1}+\frac{x}{\alpha_1 \cdots
\alpha_{n-1}}) \right| \leq (1-\frac{1}{q-1}) q^{n/2} +
\frac{1}{q-1}.\] \end{cor}

\begin{proof} The paper \cite{Ki1} proves (see the remark on page 303)
by very elementary means (using only the Bruhat decomposition of
$SL(n,q)$ with respect to a parabolic subgroup) that

\[ q^{{n \choose 2}} \sum_{\alpha_1,\cdots,\alpha_{n-1} \in F_q^*}
\lambda(\alpha_1+\cdots+\alpha_{n-1}+\frac{x}{\alpha_1 \cdots
\alpha_{n-1}}) = \sum_{g \in GL(n,q) \atop det(g)=x} \lambda(tr(g)).\]
Clearly

\[ \sum_{g \in GL(n,q) \atop det(g)=x} \lambda(tr(g)) = \frac{1}{q-1}
\sum_{\chi} \bar{\chi}(x) \sum_{g \in GL(n,q)} \chi(det(g))
\lambda(tr(g))\] where the sum is over all multiplicative characters
of $F_q^*$. From Theorem \ref{expGL} and the equation
$|G(\chi,\lambda)|=q^{1/2}$ for non-trivial $\chi,\lambda$ it follows
that

\[ \frac{1}{q^{n \choose 2}} \left|\frac{1}{q-1} \sum_{\chi}
\bar{\chi}(x)\sum_{g \in GL(n,q)} \chi(deg(g)) \lambda(tr(g)) \right|
\leq \frac{q-2}{q-1} q^{n/2} + \frac{1}{q-1}.\] Note that the second
term on the right-hand side arises from $\chi$ trivial. \end{proof}

{\bf Remark:} We observe that Fourier analysis combined with Theorem
\ref{expGL} gives a proof of Kim's relation \[ q^{{n \choose 2}}
\sum_{\alpha_1,\cdots,\alpha_{n-1} \in F_q^*}
\lambda(\alpha_1+\cdots+\alpha_{n-1}+\frac{x}{\alpha_1 \cdots
\alpha_{n-1}}) = \sum_{g \in GL(n,q) \atop det(g)=x} \lambda(tr(g))\]
which avoids the Bruhat decomposition of $SL(n,q)$ with respect to a
parabolic subgroup. Indeed, Fourier analysis gives that the left-hand
side is equal to \[ \frac{1}{q-1} q^{{n \choose 2}} \sum_{\chi}
\bar{\chi}(x) G(\lambda,\chi)^n.\] Theorem \ref{expGL} implies that this
is equal to \[ \sum_{g \in GL(n,q) \atop det(g)=x} \lambda(tr(g)).\]

	To close this section we find a generating function for \[
\sum_{g \in GL(n,q)} \lambda(x tr(g) + y tr(g^{-1})) \] with $x,y \neq
0$. The sum was evaluated by Kim (page 64 of \cite{Ki4}) using
induction and the Bruhat decomposition.

	Let $K_{\lambda}(x,y)$ with $x,y \neq 0$ denote the
Kloosterman sum $\sum_{\alpha \in F_q^*} \lambda(x \alpha +
\frac{y}{\alpha})$. Define $\tau(P)= \lambda(x tr(P) + y tr(P^{-1}))$
where $tr(P)$ is the sum of the roots of $P$ and $tr(P^{-1})$ the sum
of the reciprocals of the roots of $P$. To be explicit, given
$P=x^n+c_{n-1}x^{n-1}+ \cdots +c_0$, define $\tau(P)$ to be
$\lambda(-x c_{n-1}-y c_1/c_0)$ for $deg(P) \geq 2$ and to be
$\lambda(-x c_0-y/c_0)$ for $deg(P)=1$.  Note that $\tau(P_1P_2) =
\tau(P_1) \tau(P_2)$.

\begin{lemma} \label{Lf} For $x,y \neq 0$, \[ \prod_{\phi \neq z}
\left(\frac{1}{1- \tau(\phi) u^{deg(\phi)}/q^{i \cdot deg(\phi)}}
\right) = 1+\frac{u K_{\lambda}(x,y)}{q^i}+\frac{q u^2}{q^{2i}}.\]
\end{lemma}

\begin{proof} The left hand side is equal to \[ \sum_{P:P(0) \neq 0}
\frac{u^{deg(P)} \tau(P)}{q^{i \cdot deg(\phi)}}.\] The terms
corresponding to $deg(P) \geq 3$ all vanish. The degree 1 term is
equal to $\frac{u K_{\lambda}(x,y)}{q^i}$. The degree 2 term is equal
to

\[ \frac{u^2}{q^{2i}} \sum_{c_0 \in F_q^*, c_1 \in F_q} \lambda(-xc_1)
\lambda(-y c_1/c_0) \] which simplifies to $\frac{q
u^2}{q^{2i}}$. \end{proof}

	Given Lemma \ref{Lf}, Theorem \ref{KloostGL} is proved exactly
as Theorem \ref{expGL}.

\begin{theorem} \label{KloostGL} For $\lambda$ a non-trivial additive
character of $F_q$, and $x,y$ nonzero elements of $F_q$,

\[ 1+\sum_{n=1}^{\infty} \frac{u^n}{|GL(n,q)|} \sum_{g \in GL(n,q)}
\lambda(x tr(g) + y tr(g^{-1})) = \prod_{i \geq 1} \left( 1+\frac{u
K_{\lambda}(x,y)}{q^i} + \frac{qu^2}{q^{2i}} \right). \] \end{theorem}

\section{Gauss Sums for $U(n,q)$} \label{GaussU}

	This section derives a generating function for \[ \sum_{g \in
U(n,q) \subset GL(n,q^2)} \chi(det (g)) \lambda( tr (g)) \] where
$\chi$ is a multiplicative character of $F_{q^2}^*$ and $\lambda$ is a
nontrivial additive character of $F_{q^2}$. It is convenient to set
$\chi(0)=0$.

	We use the notation that \[ G_1(\chi,\lambda) = \sum_{\alpha
\in F_{q^2} \atop \alpha^{q+1}=1} \chi(-\alpha) \lambda(-\alpha) \]
and that \[ G_2(\chi,\lambda) = \sum_{\alpha \in F_{q^2} \atop
\alpha^{q+1}=1} \chi(\alpha) \sum_{\beta \in F_{q^2} :
\beta^{q-1}=\alpha^q \atop or \ \beta=0} \lambda(-\beta).\] (This
notation differs from the notation in the first version of the paper
but the results are the same).

	We also use an involution which maps a polynomial \[ \phi(z) =
z^{m} + \alpha_{m-1} z^{m-1} + \cdots + \alpha_1 z + \alpha_0 \] with
$\alpha_0 \neq 0$ to \[ \tilde{\phi}(z) = z^{m} +
(\alpha_1/\alpha_0)^q z^{m-1} + (\alpha_2/\alpha_0)^q z^{m-2} + \cdots
+ (\alpha_{m-1}/\alpha_0)^q z + (1/\alpha_0)^q. \] Note that
$\tilde{\phi_1 \phi_2} = \tilde{\phi_1} \tilde{\phi_2}$. The total
number of monic degree $m$ polynomials with coefficients in $F_{q^2}$
invariant under $\tilde{}$ is $q^m+q^{m-1}$ \cite{W}. To see this note
that for $m$ odd the coefficients $\alpha_0,\cdots,\alpha_{(m-1)/2}$
determine an invariant polyonomial and that $\alpha_0$ must satisfy
$\alpha_0^{q+1}=1$ but that $\alpha_1,\cdots,\alpha_{(m-1)/2}$ can be
any elements of $F_{q^2}$. The case of $m$ even is similar. All
irreducible $\phi$ invariant under this involution have odd degree
(see \cite{F1} or \cite{NP} for a proof) and the number of such
polynomials is computed in \cite{F1}.
	
	The paper \cite{F1} develops analogs of the cycle index for
the finite classical groups, based on Wall's work on its conjugacy
classes \cite{W}. For the case of the unitary groups, the conjugacy
classes correspond to the following combinatorial data. As was the
case with $GL(n,q^2)$, an element $g$ in $U(n,q)$ associaties to each
monic, non-constant, irreducible polynomial $\phi$ over $F_{q^2}$ a
partition $\nu_{\phi}$ of some non-negative integer
$|\nu_{\phi}|$ by means of the rational canonical form. This data
represents a conjugacy class if and only if $(1) \ |\nu_z|=0$,
$(2) \ \nu_{\phi}=\nu_{\tilde{\phi}}$, and $(3) \ \sum_{\phi}
|\nu_{\phi}| deg(\phi)=n$.

	Note that in the statement of Lemma \ref{L2} the second
product is over (unordered) pairs of distinct monic irreducible
polynomials which map to each other under the involution $\tilde{}$.

\begin{lemma} \label{L2} Let $\chi$ be a multiplicative character of
$F_{q^2}^*$ and let $\lambda$ be a non-trivial additive character of
$F_{q^2}$. Then

\begin{enumerate}

\item

\begin{eqnarray*}
& &  \prod_{\phi \neq z \atop \phi=\tilde{\phi}}
\left(\frac{1}{1-\chi(\phi) \lambda(\phi) u^{deg(\phi)}/q^{i \cdot
deg(\phi)}} \right) \prod_{\{\phi,\tilde{\phi}\} \atop \phi \neq
\tilde{\phi}} \left(\frac{1}{1-\chi(\phi) \chi(\tilde{\phi})
\lambda(\phi) \lambda(\tilde{\phi}) u^{2deg(\phi)}/q^{2i \cdot
deg(\phi)}} \right)\\
& = &1+u G_1(\chi,\lambda)/q^i+u^2
G_2(\chi,\lambda)/q^{2i}.
\end{eqnarray*}

\item

\begin{eqnarray*}
& & \prod_{\phi \neq z \atop \phi=\tilde{\phi}}
\left(\frac{1}{1+\chi(\phi) \lambda(\phi) u^{deg(\phi)}/q^{i \cdot
deg(\phi)}} \right) \prod_{\{\phi,\tilde{\phi}\} \atop \phi \neq
\tilde{\phi}} \left(\frac{1}{1-\chi(\phi) \chi(\tilde{\phi})
\lambda(\phi) \lambda(\tilde{\phi}) u^{2deg(\phi)}/q^{2i \cdot
deg(\phi)}} \right)\\
& = &1-u G_1(\chi,\lambda)/q^i+u^2
G_2(\chi,\lambda)/q^{2i}.
\end{eqnarray*}
\end{enumerate}
\end{lemma} 

\begin{proof} Letting $P= z^{m} + \alpha_{m-1} z^{m-1} + \cdots +
\alpha_1 z + \alpha_0$ denote a monic polynomial with coefficients in
$F_{q^2}$, the left hand side is equal to

\[ \sum_{P: \tilde{P}=P, P(0) \neq 0} \frac{\chi(P) \lambda(P)
u^{deg(P)}}{q^{i \cdot deg(P)}}.\] Observe that for $d > 2$, the
expression \[ \sum_{P: \tilde{P}=P, P(0) \neq 0 \atop deg(P)=d}
\frac{\chi(P) \lambda(P) u^{d}}{q^{i \cdot deg(P)}} \] vanishes
because (from the explicit description of invariant P) $\alpha_{d-1}$
is equidistributed over all elements of $F_{q^2}$ given the value of
$\alpha_0$. The computations for $deg(P)=1$ and $deg(P)=2$ are
straightforward.

	The second assertion follows from the first assertion by
replacing $u$ by $-u$ and using the fact that all irreducible
polynomials invariant under $\tilde{}$ have odd degree. \end{proof}

\begin{theorem} \label{produ} Let $\chi$ be a multiplicative character
of $F_{q^2}^*$ and let $\lambda$ be a non-trivial additive character of
$F_{q^2}$. Then

\begin{eqnarray*}
& & 1+\sum_{n=1}^{\infty} \frac{u^n}{|U(n,q)|} \sum_{g \in U(n,q)
\subset GL(n,q^2)} \chi(det (g)) \lambda( tr (g))\\
& = & \prod_{i \geq 1}
\left( 1+(-1)^{i+1} u G_1(\chi,\lambda)/q^i+u^2
G_2(\chi,\lambda)/q^{2i} \right).
\end{eqnarray*}

\end{theorem}

\begin{proof} Arguing as in \cite{F1} and using the fact that Gauss
sums are multiplicative on polynomials, it follows that

\[ 1+\sum_{n=1}^{\infty} \frac{u^n}{|U(n,q)|} \sum_{g \in U(n,q)}
\chi(det (g)) \lambda( tr (g)) \] is equal to

\[  \prod_{i \geq 1} \prod_{\phi \neq z \atop \phi=\tilde{\phi}}
\left(\frac{1}{1- (-1)^{i+1} \chi(\phi)
\lambda(\phi) u^{deg(\phi)}/q^{i \cdot deg(\phi)}} \right)
\prod_{\{\phi,\tilde{\phi}\} \atop \phi \neq \tilde{\phi}}
\left(\frac{1}{1-\chi(\phi) \chi(\tilde{\phi}) \lambda(\phi)
\lambda(\tilde{\phi}) u^{2deg(\phi)}/q^{2i \cdot deg(\phi)}} \right)
.\] This product can be broken down into terms according to whether
$i$ is even or odd and the result follows from Lemma
\ref{L2}. \end{proof}

	Note that if one factors the right hand side of the generating
function of Theorem \ref{produ} into linear factors in $u$ (which is
certainly possible in odd characteristic), then the right hand side
can be expanded using Euler's identity

\[ \prod_{i \geq 1} (1-\frac{u}{q^i}) = \sum_{n=0}^{\infty}
\frac{(-u)^n}{(q^n-1) \cdots (q-1)}.\]

\section{Acknowledgements} This research was supported by an NSF
Postdoctoral Fellowship. The author thanks Robin Champan for helpful
correspondence.

\end{document}